\documentclass{amsart}
\usepackage{amsfonts,amssymb,amscd,amsmath,enumerate,verbatim,calc}
\usepackage[all]{xy}

\newcommand{\CM}{Cohen-Macaulay}

\newcommand{\wrt}{with respect to}

\newcommand{\m}{\mathfrak{m} }

\newcommand{\fp}{\mathfrak{p}}

\newcommand{\QQ}{\mathbb{Q} }
\newcommand{\W}{\mathbf{W}_{\bullet} }
\newcommand{\D}{\mathbf{D}_{\bullet} }
\newcommand{\Dc}{\mathbf{D}^{\bullet} }
\newcommand{\I}{\mathbf{I^\bullet} }
\newcommand{\K}{\mathbb{K}_{\bullet} }
\newcommand{\Kc}{\mathbb{K}^{\bullet} }

\newcommand{\C}{\mathbf{C} }
\newcommand{\T}{\mathbf{T}^{\bullet}  }
\newcommand{\X}{\mathbf{X} }
\newcommand{\Y}{\mathbf{Y} }
\newcommand{\Z}{\mathbf{Z}^{\bullet}  }
\newcommand{\F}{\mathbf{F}_{\bullet} }

\newcommand{\rt}{\rightarrow}
\newcommand{\xar}{\longrightarrow}
\newcommand{\ov}{\overline}

\newcommand{\bx}{\mathbf{x}}
\newcommand{\bu}{\mathbf{u}}
\newcommand{\by}{\mathbf{y}}

\newcommand{\grade}{\operatorname{grade}}
\newcommand{\depth}{\operatorname{depth}}
\newcommand{\Tot}{\operatorname{Tot}}

\newcommand{\height}{\operatorname{height}}

\newcommand{\projdim}{\operatorname{projdim}}
\newcommand{\injdim}{\operatorname{injdim}}
\newcommand{\Hom}{\operatorname{Hom}}
\newcommand{\Ext}{\operatorname{Ext}}
\newcommand{\Tor}{\operatorname{Tor}}

\theoremstyle{plain}

\newtheorem{theorem}{Theorem}[section]
\newtheorem{corollary}[theorem]{Corollary}
\newtheorem{lemma}[theorem]{Lemma}
\newtheorem{proposition}[theorem]{Proposition}

\theoremstyle{definition}

\newtheorem{s}[theorem]{}
\newtheorem{remark}[theorem]{Remark}
\newtheorem{example}[theorem]{Example}

\theoremstyle{remark}

\begin{document}

\title{Properties of  Koszul homology modules}
 \author[Uwe~Nagel and Tony~J.~Puthenpurakal]{Uwe~Nagel${}^*$ and Tony~J.~Puthenpurakal${}^\dagger$  }
%\date{\today}
\address{ Department of Mathematics, University of Kentucky, 715 Patterson Office Tower, Lexington, KY 40506-0027, USA}
\email{uwenagel@ms.uky.edu}
\address{Department of Mathematics, Indian Institute of Technology Bombay, Powai, Mumbai 400 076, India}

\email{tputhen@math.iitb.ac.in}

\thanks{${}^*$  The work for this paper was done while the first
author was sponsored by the National Security Agency under Grant
Number H98230-07-1-0065. \\
${}^\dagger$ The work for this paper was done while the second
author was visiting University of Kentucky by a fellowship from Department of Science and Technology, India
}

\begin{abstract}
We investigate various module-theoretic properties of Koszul homology under mild conditions. These include their depth, $S_2$-property and their Bass numbers
\end{abstract}
\maketitle

\section{Introduction}

This note is an attempt to study module theoretic properties of Koszul homology modules.
Let $R$ be a commutative Noetherian  ring and let $M$ be a finitely generated $R$-module. Let $I$ an $R$-ideal in $R$  and let $\by = \{y_1,\ldots,y_l \}$ be any  set of  generators of $I$ and assume that $M \neq IM$. Denote by $\K = \K(\by, M)$  the Koszul complex associated to $\by$ with coefficients in $M$. Set $H_i(\by, M)$ the $i^{th}$ Koszul homology module of $M$ \wrt\ $\by$. It is well-known that the Koszul complex is grade sensitive, that is, if $\grade(I, M) = g$ then $H_{l-g} (\by, M) \neq 0$ and $H_i (\by, M) = 0$ whenever $i > l-g$.
Furthermore, if $(R, \m)$ is a local \CM\ ring and $I$ is an  ideal, then  all  non-zero
 $H_i(\by, R)$ have the same dimension as $A = R/I$ (see, e.g., \cite[4.2.2]{Vas}). However, apart from special cases (see, e.g., \cite{Green}) not much is known about Koszul homology modules. Indeed, W. Vasconcelos writes in his book ``Integral closures'' \cite[page 280]{VasBook-05}:
\begin{quote}
"While the vanishing of the homology of a Koszul complex $K(\bx, M)$ is easy to track, the module theoretic properties
of its homology, with the exception of the ends, is difficult to fathom. For instance, just trying to see
whether a prime is associated to some  $H_i(\bx, M)$  can be very hard."
\end{quote}
The purpose of this note is to enhance our knowledge about Koszul homology by establishing the following five results:

We first give the following general estimate for depth of $H_{l-g-1}(\by, R)$.
\begin{theorem}\label{depth-estimate}
Let $(R,\m)$ be a local \CM \ ring and let $I$ be an  ideal. Set $g = \grade(I)$. Then
$$ \depth H_{l-g-1}(I) \geq \min \{ 2, \depth H_{l-g}(I) - 2 \}. $$
\end{theorem}

It is well-known that if $g = \grade(I,M)$ then   $H_{l-g}(\by,M)
\cong  \Ext^{g}_{R}(R/I, M)$. We show:

\begin{theorem}\label{main-perfect}
Let $(R,\m)$ be a \CM \  local ring and let $I$ be a perfect ideal.
Let $M$ be a  maximal \CM \ $A$-module. Set $g = \grade(I)$. Then
 $$ H_{l-g-1}(\by, M) \cong \Ext^{g}_{R}(H_1(\by), M). $$
 In particular $ H_{l-g-1}(\by, M)$ satisfies $S_2$.
\end{theorem}
If $(R,\m)$ is \CM \ and $I$ is a \CM \ ideal,  Vasconcelos notes
that $H_{n-g-1}(\by, R)$ is $S_2$; see \cite[1.3.2]{VasKH}.  We
prove:

\begin{theorem}\label{main-Gor}
 Let $(R,\m)$ be a  Gorenstein and $I$ be a \CM \ ideal.  Then
 \[
 H_{n-g-1}(\by) \cong \Ext_R^g(H_1(\by), R).
 \]
 \end{theorem}

Next we consider the projective dimension of Koszul homology
modules. Note that if $I$ is perfect then $H_{l-g}(\by,R) \cong
\Ext^{g}_{R}(R/I, R)$ has finite projective dimension.
\begin{theorem}\label{proj-finite}
Let $(R,\m)$ be a local Noetherian and let $I$ be a   perfect ideal
with $\nu(I) = \grade(I) + 2 $, where $\nu (I)$ denotes the number
of minimal generators of $I$.  Then $\projdim_R H_1(I)$ is finite.
\end{theorem}

We conclude with an estimate of certain Bass numbers of $H_1(I)$. Let $\nu(E)$ denote the  minimal number of generators of an $R$-module $E$ and let $\mu_i(\m, E) = \ell\left( \Ext^{i}_{R}(k,E) \right)$ denote the $i$-th Bass number of $E$ (\wrt \ $\m$). If $(R,\m)$ is regular local of dimension
$d$ then one can verify $\mu_d(\m, R/I) = 1$.
\begin{theorem}\label{bass-estimate}
Let $(R,\m)$ be a regular local ring and let $I$ be an ideal in $R$. Then
$$ l - \mu_{d-1}(\m, R/I) \leq \mu_d(\m, H_1(I)) \leq l - \mu_{d-1}(\m, R/I)  + \mu_{d-2}(\m, R/I).$$
\end{theorem}

We now briefly describe  the organization of this paper. In Section
2 we introduce some notation. In Section 3 we prove Theorems
\ref{main-perfect}, \ref{main-Gor} and \ref{proj-finite}. In Section
4 we prove Theorem \ref{depth-estimate}. We conclude with an
estimate of certain Bass numbers in Section 5.

\section{Notation}
In this paper all rings are commutative Noetherian.   Let $R$ be a
ring, $I$   an ideal in $R$ and let $M$ be  an  $R$-module
(not-necessarily finitely generated).

\s Let $\by = \{y_1,\ldots,y_l \}$ be a  set of  generators of $I$ and let $\K(\by, R)$ be the Koszul complex with respect  to $\by$. Set
 $$\K = \K(\by, R)\colon \quad 0\rt K_l \rt \cdots K_1 \rt K_0 \rt 0. $$
 Let $\Kc = \Hom_R(\K,R)$ be the Koszul co-chain complex \wrt\ $\by$. Let $\K(\by,M)$ and  $\Kc(\by, M) $ be  respectively the Koszul complex  and co-chain complex  \wrt \ $\by$ with coefficients in $M$.

\s If $\D$ is a chain complex of $R$-modules then we set $H_{i}(\D)$ to denote its $i^{th}$ homology module. Likewise if $\Dc$ is a co-chain complex of $R$-modules then we set $H^{i}(\Dc)$ to denote its $i^{th}$ cohomology module. Let $H_{i}(\by, M)$ and $H^{i}(\by, M)$ denote respectively the $i^{th}$  Koszul homology and cohomology module \wrt \ $\by$ with coefficients in  $M$.

\section{Proof of  theorems \ref{main-perfect}, \ref{main-Gor} and \ref{proj-finite}}
 In this section we  prove Theorems \ref{main-perfect}, \ref{main-Gor} and \ref{proj-finite}

\s\label{setup-doublecomplex}
  Let $\I$ be a "deleted" injective  resolution of $M$;
$$ \I \colon \quad 0 \rt I^0 \rt I^1 \rt \cdots \rt I^n \rt I^{n+1} \rt \cdots $$
We consider the Hom co-chain complex $\C = \Hom(\K, \I)$; see \cite[2.7.4]{Weibel}. Set
\[
\C = \{C^{pq} \}_{p,q\geq 0} \quad \text{where} \quad  C^{pq} = \Hom(K_p,I^q).
\]
Set $\T  = \Tot^{\bigoplus}(\C)$ where
$$ T^{n} = \bigoplus_{p+q = n}C^{pq}.$$

\begin{proposition}\label{I-seq}[with hypothesis as in \ref{setup-doublecomplex}]
The spectral sequence $\{ {}^{I}E_{r}^{pq} \}$ collapses; hence,  for each $i \geq 0$, we have
$H^{i}(\T) = H^{i}(\by, M)$.
\end{proposition}
\begin{proof}
Notice
\begin{align*}
{}^{I}E_{1}^{pq} &= H^q(\Hom_R(K_p, \I) \\
                 &= \Ext^{q}_{R}(K_p, M) \\
                 &= \left\{
                      \begin{array}{ll}
                       \Hom_{R}(K_p, M) , & \text{for} \  q = 0 ; \\
                        0 & \text{for} \ q > 0.
                      \end{array}
                    \right.
\end{align*}
The last equality is  true since $K_p$ is a free  $R$-module. It follows that
\[
{}^{I}E_{2}^{pq} = \left\{
                      \begin{array}{ll}
                       H^p(\by, M) , & \text{for} \  q = 0 ; \\
                        0 & \text{for} \ q > 0.
                      \end{array}
                    \right.
\]
Hence the spectral sequence collapses at $E_2$ and the claim follows.
\end{proof}
\begin{proposition}\label{II-seq}[with hypothesis as in \ref{setup-doublecomplex}]
$${}^{II}E_{2}^{pq}  = \Ext^{p}_{R}(H_q(\by), M).$$
\end{proposition}
\begin{proof}
Notice
\begin{align*}
{}^{II}E_{1}^{pq}  &= H^{q}\left(\Hom_R(\K, I^p) \right) \\
 &= \Hom_R\left(H_{q}(\by), I^p)\right).
\end{align*}
The last equality is true since $I^p$ is an injective $R$-module. Thus
\begin{align*}
{}^{II}E_{2}^{pq}  &= H^{p}\left(\Hom_R(H_{q}(\by), \I ) \right) \\
 &= \Ext^{p}_{R}\left(H_{q}(\by), M)\right);
\end{align*}
as claimed.
\end{proof}

\s \label{cases}\textbf{ ${}^{II}E_{2}$-page in three Special Cases:}
We will consider the following three special cases

\begin{enumerate}[(i)]
\item
\label{case-Gor}   \emph{$(R,\m)$ is a $d$-dimensional Gorenstein local ring, $A = R/I$ is CM and $M = R$.  }
\begin{enumerate}[\rm (a)]
  \item
By local duality we get $\Ext^{p}_{R}(H_{0}(I), M) = 0$ for $p>g$. So
 $${}^{II}E_{2}^{p,0}  = 0 \quad \text{for} \ p > g.$$
 \item
 Since $\grade H_q(\by) = g$ for all $q$ and  $M = R$ is maximal \CM.  So
  $${}^{II}E_{2}^{pq}  = 0 \quad \text{for} \ p < g.  $$
\end{enumerate}
\item
\label{case-perfect}  \emph{$(R,\m)$ is a $d$-dimensional \CM \ local ring with a canonical module, the ideal $I$ is perfect (in particular $A = R/I$ is CM) and $M $ is a Maximal \CM \ $R$-module.  }
\begin{enumerate}[\rm (a)]
  \item
Since $H_0(\by) = R/I$ has projective dimension $g$  we get that \\ $\Ext^{p}_{R}(H_{0}(I), M) = 0$ for $p>g$. So
 $${}^{II}E_{2}^{p,0}  = 0 \quad \text{for} \ p > g.$$
 \item
 Since $\grade H_q(\by) = g$ for all $q$ and  $M$  is maximal \CM \ we get
  $${}^{II}E_{2}^{pq}  = 0 \quad \text{for} \ p < g.  $$
\end{enumerate}
\item\label{G2}
\emph{$(R,\m)$ is Noetherian local of dimension $d$ and $M$ is a finitely generated $R$-module. The ideal $I$ is perfect and has
$l= \mu(I) = \grade(I) +2$. Finally, $\by$ is a minimal set of generators of $I$.}\\
\begin{enumerate}
  \item
  We have that $H_{2}(I) \cong \Ext^{g}_{R}(R/I, R)$ has projective dimension $g$. So $\Ext^{p}_{R}(H_{2}(I), M) = 0$ for $p>g$.
  $${}^{II}E_{2}^{p,l-g}  = 0 \quad \text{for} \ p > g.$$
  \item
  Similarly as $R/I$ has projective dimension $g$ we get $\Ext^{p}_{R}(H_{0}(I), M) = 0$ for $p>g$.  So
  $${}^{II}E_{2}^{pq}  = 0 \quad \text{for} \ p < g.  $$
\end{enumerate}
\end{enumerate}

We now establish Theorem \ref{main-perfect}.
\begin{proof}[Proof of Theorem \ref{main-perfect}]
Assume  $R,I,A, M$ are as in \ref{cases}(\ref{case-perfect}).
Recall ${}^{II}E_{r}$ has differential of degree $(r,-r+1)$.
Using the vanishing results in \ref{cases}(\ref{case-perfect}) we get that
\begin{align*}
{}^{II}E_{\infty}^{g,1} &= {}^{II}E_{2}^{g,1} = \Ext^{g}_{R}(H_{1}(\by), M) \quad \text{and} \\
{}^{II}E_{\infty}^{g+1,0} &= {}^{II}E_{2}^{g+1,0} = 0.
\end{align*}
The only non-zero term in   ${}^{II}E_{\infty}$ with total degree $ g+1$ is ${}^{II}E_{\infty}^{g,1}$.
It follows that
\[
H^{g+1}(\T) \cong {}^{II}E_{\infty}^{g,1} =  \Ext^{g}_{R}(H_{1}(\by), M) .
\]
Proposition \ref{I-seq} provides our claim.
In particular $ H_{l-g-1}(\by, M)$ is $S_2$.
\end{proof}

We now prove   our second main result.
\begin{proof}[Proof of Theorem\ref{main-Gor}]
Assume  $R,I,A$ are as in  \ref{cases}(\ref{case-Gor}).
We use the spectral sequence with $M = R$.
Recall ${}^{II}E_{r}$ has differential of degree $(r,-r+1)$.
Using the vanishing results in \ref{cases}(\ref{case-Gor}) we get that
\begin{align*}
{}^{II}E_{\infty}^{g,1} &= {}^{II}E_{2}^{g,1} = \Ext^{g}_{R}(H_{1}(\by), R) \quad \text{and} \\
{}^{II}E_{\infty}^{g+1,0} &= {}^{II}E_{2}^{g+1,0} = 0.
\end{align*}
The only non-zero term in   ${}^{II}E_{\infty}$ with total degree $ g+1$ is ${}^{II}E_{\infty}^{g,1}$.
It follows that
\[
H^{g+1}(\T) \cong {}^{II}E_{\infty}^{g,1} =  \Ext^{g}_{R}(H_{1}(\by), R) .
\]
Proposition \ref{I-seq} provides our claim.
\end{proof}

We
now establish Theorem \ref{proj-finite}
\begin{proof}[Proof of Theorem \ref{proj-finite}]
Recall ${}^{II}E_{r}$ has differential of degree $(r,-r+1)$.
Also recall that $g \leq d$.

Using the vanishing results in \ref{cases}{iii} we get
\[
{}^{II}E_{\infty}^{d+3,1} = {}^{II}E_{2}^{d+3,1} = \Ext^{d+3}_{R}(H_{1}(\by), M).
\]
By Proposition \ref{I-seq};  $H^{d+4}(\T) = 0$. It follows that
\[
\Ext^{d+3}_{R}(H_{1}(\by), M) = 0,
\]
where the module $M$ is an arbitrary $R$-module. So we get $\projdim_R H_{1}(I) < \infty$.
\end{proof}

An easy consequence to \ref{proj-finite} is the following result.
\begin{corollary}\label{nice-app}
Let $(R,\m)$ be a Gorenstein local ring and let $I$ be a perfect ideal in $R$ with $\nu(I) = \grade(I) + 2$. Set $A = R/I$. Then
$H_1(I)$ is a perfect $R$-module and a self-dual $A$-module
 \end{corollary}
\begin{proof}
By a result due to Avramov and Herzog  \cite[Supplement]{AH}, the ideal $I$ is strongly CM.
Using \ref{proj-finite} we get that $H_1(I)$ is a perfect $R$-module.

Let $\omega$ be the canonical module of $A$.
By \ref{main-Gor}
we have
$H_1(I) \cong \Ext^{g}_R(H_1(I),R)$. Notice
$$\Ext^{g}_R(H_1(I),R) \cong \Hom_A (H_1(I), \omega).$$
Thus $H_1(I)$ is a self-dual $A$-module.
\end{proof}
\section{Proof of Theorem \ref{depth-estimate}}
Throughout this section
$(R,\m)$ is \CM \ local ring of dimension $d$ and $I$ is an ideal in $R$. Let $\bu = u_1,\ldots,u_l $ be a system of minimal generators for $I$.
Let $\K = \K(\bu,R)$ be the Koszul complex \wrt \ $\bu$.
Let $g = \grade(I)$ and let $\bx = x_1,\ldots,x_{d-g} \in \m$ be such that
\begin{enumerate}
  \item $\bx$ is a $R$-regular sequence.
  \item $\bx$ is a system of parameters for $R/I$.
\end{enumerate}

Let $\C^\bullet$ be the \v{C}ech complex on $\bx$. We write $\C^\bullet[-(d-g)]$ homologically and call it $\D$.
So
\[
\D \colon 0\rt D_{d-g} \rt \cdots \rt D_1 \rt D_0 \rt 0
\]
and $H_i(D\otimes M) = H^{d-g-i}_{\bx}(M)$
for a $R$-module $M$.

Consider the double complex $\X = \K \otimes \D$ and set $\W = \Tot(\X)$.
We look at the two standard spectral sequences associated to $\X$.

\begin{proposition}\label{I-seq-depth}
The spectral sequence $\{ {}^{I}E^{r}_{pq} \}$ collapses; hence,  for each $i \geq 0$, we have
$H_{i}(\W) = H_{i}(\bu, H^{d-g}_{\bx}(R))$.
\end{proposition}
\begin{proof}
$ {}^{I}E^{0}_{pq} = K_p \otimes D_q $. So
\begin{align*}
 {}^{I}E^{1}_{pq} &= H_q(K_p \otimes \D) \\
             &= H^{d-g-q}_{\bx}(K_p) \\
              &= H^{d-g-q}_{\bx}(R) \otimes K_p \\
              &= \begin{cases} 0 & \text{for} \ q > 0; \\    H^{d-g}_\bx(R)\otimes K_p &   \text{for} \ q = 0.         \end{cases}
\end{align*}
Therefore
\[
 {}^{I}E^{2}_{pq} = \begin{cases} 0 & \text{for} \ q > 0; \\    H_p(\bu, H^{d-g}_\bx(R)) &   \text{for} \ q = 0.         \end{cases}
\]
The result follows.
\end{proof}

\begin{proposition}\label{II-seq-depth}
$${}^{II}E^{2}_{pq}  = H^{d-g-p}_{\m}(H_q(I)).$$
\end{proposition}
\begin{proof}
$ {}^{II}E^{0}_{pq} = K_q \otimes D_p $. So
\begin{align*}
 {}^{II}E^{1}_{pq} &= H_q(\K \otimes D_p) \\
                   &= H_q(\K) \otimes D_p; \quad \text{since $D_p$ is a flat $R$-module} \\
                   &= H_q(I) \otimes D_p
\end{align*}
Therefore
\begin{align*}
 {}^{II}E^{2}_{pq} &= H_p(H_q(I) \otimes \D) \\
      &=  H^{d-g-p}_{\bx}(H_q(I)) \\
         &=  H^{d-g-p}_{\m}(H_q(I)).
\end{align*}
\end{proof}
Surprisingly we have the following vanishing result.
\begin{proposition}\label{vanish-depth}Adopt the above assumptions. Then
$H_i(I, H^{d-g}_\bx(R)) = 0$ for $i > l -g$.
\end{proposition}
To prove this result the following Lemma is needed.
\begin{lemma}\label{weak}
Let $(R,\m)$ be a \CM \ local ring. Let
$x_1,\ldots,x_r, y_1,\ldots,y_s$ be an $R$-regular sequence. Then
\begin{enumerate}[\rm (1)]
  \item $\by = y_1,\ldots,y_s$ is a weak
$H^r_\bx(R)$-regular sequence.
  \item $H^r_\bx(R)/\by H^r_\bx(R) = H^r_{\bx}(R/\by R)$
\end{enumerate}
\end{lemma}
\begin{proof}
It is sufficient to prove it for $s = 1$. Set $y = y_1$ and $\ov R = R/yR$. Consider the
exact sequence
\[
0 \xar R \xrightarrow{y} R \xar \ov R \xar 0
\]
Notice $\bx$ is a $R\oplus \ov R$-regular sequence. Therefore taking local cohomology \wrt \ $\bx$ we obtain
\[
0 \xar H^r_\bx(R) \xrightarrow{y}  H^r_\bx(R) \xar  H^r_\bx(\ov R) \xar 0.
\]
So $y$ is $H^r_\bx(R)$-regular and $H^r_\bx(R)/yH^r_\bx(R) = H^r_\bx(\ov R)$.
\end{proof}

\begin{proof}[Proof of Proposition \ref{vanish-depth}]
Choose $\by = y_1,\ldots,y_g$ in $I$ such that $\bx,\by$ is a system of parameters for $R$ and hence a $R$-regular sequence.
The result follows from Lemma \ref{weak} and \cite[1.6.16]{BH}.
\end{proof}

The following example shows that the result in Theorem \ref{depth-estimate}  cannot be improved in general.
\begin{example}
Let $R = \QQ[x,y,z,w, a,b,c,d]$ and let $I $ be the ideal generated by the maximal minors of
$\psi$  where
\[
\psi= \begin{pmatrix} a & b & c & d \\ x & y & z & w \end{pmatrix}.
\]
By  \cite{Eagon-N},  we get that
$\height I = 3$, so $\dim A = 5$. Let $\by $ be the set of  minimal generators of $I$. So $ l=6, g= 3$.
Using MACAULAY \cite{M2}, one  verifies that $\depth H_{2}(\by) = 2$. Here $\depth H_{3}(\by) = 5$.
\end{example}

We now give the proof of Theorem \ref{depth-estimate}.
\begin{proof}
We have to show the following
\begin{enumerate}[\rm(1)]
  \item If $\depth H_{l-g}(I)\geq 3$ then $\depth H_{l-g-1}(I) \geq 1$.
  \item  If $\dim H_{l-g}(I)  \geq 4$ then $\depth H_{l-g-1}(I) \geq 2$.
\end{enumerate}
We use the two standard spectral sequences induced on the above double complex $\X = \K \otimes \D$. Recall ${}^{II}E^{r}$ has differential of degree $(-r,r-1)$. We also notice that
$$ {}^{II}E^{2}_{pq}  = H^{d-g-p}_{\m}(H_q(I)) = 0 \quad \text{for} \ p > d-g.$$
If  $\depth H_{l-g}(I)\geq i+1$ then
$$ {}^{II}E^{2}_{d-g-j, l-g}  = 0\quad \text{for} \ j = 0,\ldots,i. $$

(1) If  $\depth H_{l-g}(I)\geq 3$ then using the above vanishing results   we get that
\[
{}^{II}E^{\infty}_{d-g,l-g-1} = {}^{II}E^{2}_{d-g,l-g-1} =  H^{0}_{\m}(H_{l-g-1}(I))
\]
Since $\dim  H_{l-g}(I) = d-g  \geq  \depth  H_{l-g}(I) = 3$, the total degree of \\
$ {}^{II}E^{\infty}_{d-g,l-g-1}$ is $r = d-g + l-g-1 \geq l-g+2$. As
 $ {}^{II}E^{\infty}_{d-g,l-g-1}$ is a subquotient of $H_r(\W) =0$ (by Proposition \ref{vanish-depth}), we get that $ {}^{II}E^{\infty}_{d-g,l-g-1} = 0.$
 Thus
 $  H^{0}_{\m}(H_{l-g-1}(I))  = 0$.  Therefore
 $\depth H_{l-g-1}(I) \geq 1$.

 (2)   If $\depth H_{l-g}(I)\geq 4$ then similarly as above we get
\[
{}^{II}E^{\infty}_{d-g-1,l-g-1} = {}^{II}E^{2}_{d-g-1,l-g-1} =  H^{1}_{\m}(H_{l-g-1}(I))
\]
The total degree of $ {}^{II}E_{\infty}^{d-g-1,l-g-1}$ is $r = d-g-1 + l-g- 1 \geq l-g+2$. By an argument  similar to (1) it follows
that $H^{1}_{\m}(H_{l-g-1}(I))  = 0$.   By (1) we also have that $H^{0}_{\m}(H_{l-g-1}(I))  = 0$.   Therefore
 $\depth H_{l-g-1}(I) \geq 2$.
\end{proof}
\section{Bass numbers}
In this section $(R,\m)$ is a Gorenstein local ring.
Let $\nu(E)$ denote the  minimal number of generators of an $R$-module $E$ and let $\mu_i(\m, E) = \ell\left( \Ext^{i}_{R}(k,E) \right)$ denote the $i$-th Bass number of $E$ (\wrt \ $\m$).
\begin{theorem}\label{bass}
Let $(R,\m)$ be a Gorenstein local ring of dimension $d$. Set $l = \nu(I), g = \grade(I)$ and assume that $l \geq g+2$. We have the following
\begin{enumerate}[\rm (I)]
\item
 Assume $I$ be a strongly \CM \ ideal in $R$. Set $c = d-g$.
Then
$$ l - \mu_{c+1}(\m, H^{g}(I))   \leq  \mu_{c}(\m, H^{g+1}(I)) \leq l - \mu_{c+1}(\m, H^{g}(I)) + \mu_{c+2}(\m, H^{g}(I)). $$
\item
Assume $\projdim_R H_i(I)$ is finite for all $i$. (Notice $I$ need not be strongly \CM). Then
 \begin{enumerate}[\rm (a)]
 \item
 $ \mu_{d}(\m, R/I) = 1$.
 \item
 $\displaystyle{ l - \mu_{d-1}(\m, R/I) \leq \mu_d(\m, H_1(I)) \leq l - \mu_{d-1}(\m, R/I)  + \mu_{d-2}(\m, R/I).}$
 \end{enumerate}
\end{enumerate}
\end{theorem}
\begin{proof}
Let $\F$ be a "deleted" minimal free resolution of $k = R/\m$, let $\I$ be a "deleted" minimal injective resolution of $R$ and let
$\K$ be the Koszul complex on a set of minimal generators of $I$.

Consider the double co-chain complexes
\[
\X = \Hom_R (\Tot(\F\otimes_R \K), \I)\quad \text{and} \quad \Y = \Hom_R (\F,  \Tot( \Hom_R(\K, \I))).
\]
Since all complexes involved  are first quadrant complexes we have $\X \cong \Y$; cf. \cite[2.7.3]{Weibel}.
Both the cases considered involve computing the cohomology of $\Z = \Tot(\X)$.  We use the second standard spectral sequence associated to $\X$
to compute cohomology of $\Z$.
Notice  $\Tot(\X) \cong \Tot(\Y)$. We use the first standard spectral spectral sequence on $\Y$ to derive our results.

Set $\D = \Tot(\F\otimes_R \K)$ and $\T = \Tot( \Hom_R(\K, \I))$. By Proposition \ref{I-seq} we have that
$H^i(\T) = H^i(I)$; the $i$-th Koszul cohomology of $I$.

Now we compute the homology of $\D$. We use the second standard spectral sequence on $\F\otimes_R \K$.
So $E^0_{pq} = F_q \otimes K_p$. So we get $E^1_{pq} = H_q(\F\otimes K_p) = \Tor^R_q(k, K_p)$. Since
$K_p$ is free we have
\[
E^1_{pq} = \begin{cases} 0 & \text{for}\ q \neq 0 \\     k\otimes K_p    & \text{for}\ q = 0               \end{cases}
\]
It follows that
\[
E^2_{pq} = \begin{cases} 0 & \text{for}\ q \neq 0 \\     k^{\binom{l}{p}}    & \text{for}\ q = 0               \end{cases}
\]

Now we compute the  cohomology of $\Z = \Tot(\X)$,
 $$\X = \Hom_R (\Tot(\F\otimes_R \K), \I)$$
  and using the second standard spectral sequence
for $\X$. So $E^{pq}_{0} = \Hom_R(D_q, I^p)$. Therefore
\begin{align*}
E^{pq}_{1} &= H^q\left( \Hom(\D, I^p  \right) \\
           &= \Hom_R \left(H_q(\D), I^p   \right); \quad \text{since $I^p$ is injective}, \\
           &= \Hom_R \left( k^{\binom{l}{q}}, I^p   \right) \\
           &= \Hom_R \left( k , I^p   \right)^{\binom{l}{q} }
\end{align*}
Therefore
\begin{align*}
E^{pq}_{2} &= H^p\left( \Hom(k, \I \right)^{\binom{l}{q} } \\
           &= \Ext^p_R(k,R)^{\binom{l}{q} } \\
           &= \begin{cases} 0 & \text{for}\ p \neq d \\     k^{\binom{l}{q}}    & \text{for}\ p = d               \end{cases}
\end{align*}
Thus this spectral sequence collapses.
It follows that
\begin{equation}\label{total}
 H^i(\Z) = \begin{cases}      k^{\binom{l}{d-i}}    & \text{for}\ d \leq i \leq d+l \\ 0 & \text{otherwise}             \end{cases}
\end{equation}

We now use the fact that $\X \cong \Y$. So $\Z \cong \Tot(\Y)$.  We compute the cohomology of $\Z$ by using the first standard
spectral sequence on $\Y$.
So $E_0^{pq} = \Hom_R(F_p, T^q)$. Therefore we get
\begin{align*}
E^{pq}_{1} &= H^q\left( \Hom_R(F_p, \T)\right) \\
  &= \Hom_R \left( F_p, H^q(\T) \right); \quad \text{since $F_p$ is free} \\
   &= \Hom_R \left( F_p, H^q(I) \right)
\end{align*}
Therefore
\begin{align*}
E^{pq}_{2} &= H^p\left( \Hom_R \left( \F, H^q(I) \right) \right) \\
&= \Ext_{R}^{p}(k, H^q(I))
\end{align*}

Now we distinguish the two cases considered in the statement of the theorem

Case \rm{(I)}:  Assume the ideal $I$ is strongly \CM. \\
Since $H^q(I) = 0$ for $q <g$ we have $E^{pq}_{2} = 0$ for $q <g$. Also since if $H_q(I) \neq 0$ it is a \CM \ $R$-module of dimension
$d-g$ we get $E^{pq}_{2} = 0$ for $p < d-g$.

 We look at elements of total degree $ d+1$
There are only two terms of total degree $d+1$. They will make up the filtration for $H^{d+1}(\Z)= k^l$.
So
\begin{equation}\label{t}
l = \ell\left( E^{d-g+1,g}_{\infty} \right) +  \ell \left(E^{d-g,g+1}_{\infty} \right).
\end{equation}
Notice
\begin{enumerate}[\rm (i)]
  \item $E^{d-g+1,g}_{\infty} = E^{d-g+1,g}_{2} = \Ext_{R}^{d-g+1}(k, H^g(I))$.
  \item We look at $E^{pq}_{3}$ when $p = d-g$ and $q = g+1$. Notice we have an exact sequence
\begin{equation}\label{w}
0 \xar E^{d-g,g+1}_{3} \xar    E^{d-g,g+1}_{2}  \xar E^{d-g+2,g}_2
\end{equation}
Recall that $E_r$ has differential of degree $(r,-r+1)$. It follows that
\[
 E^{d-g,g+1}_{3} =  E^{d-g,g+1}_{\infty}
\]
\end{enumerate}
The result follows.

Case \rm{(II)}: Assume $\projdim_R H_i(I)$ is finite for all $i$.\\
Since $R$ is Gorenstein we have $\injdim_R H_i(I)  = d$  for all $i$. It follows that
\[
E^{pq}_{2} = \Ext_{R}^{p}(k, H^q(I)) = 0 \quad \text{for all}\ \ p > d.
\]
(a.) The only term with total degree $ d+l$ is $E^{d,l}_{\infty} $. Notice
\[
E^{d,l}_{\infty} = E^{d,l}_{2} =  \Ext_{R}^{d}(k, H^l(I))
\]
Now the result follows from \ref{total}.

(b) We look at the two terms of total degree $d+l-1$. The proof is almost similar to that of Case  \rm{(I)}. Except that here we have an exact sequence
\[
\Ext_R^{d-2}(k, H^l(I)) \xar \Ext_{R}^{d}(k, H^{l-1}(I)) \xar E^{d, l-1}_{\infty} \xar 0.
\]
\end{proof}
\begin{remark}
We wonder if it is possible to relax the assumption in (II) and still have the conclusion that $\mu_d(\m, R/I) = 1$.
\end{remark}
We now give a proof of Theorem \ref{bass-estimate}.
\begin{proof}[Proof   of Theorem \ref{bass-estimate} ]
Since $R$ is regular local we have $\projdim_R H_i(I)$ is finite for all $i$. The result follows from
Theorem \ref{bass}
\end{proof}
\providecommand{\bysame}{\leavevmode\hbox to3em{\hrulefill}\thinspace}
\providecommand{\MR}{\relax\ifhmode\unskip\space\fi MR }
\providecommand{\MRhref}[2]{%
  \href{http://www.ams.org/mathscinet-getitem?mr=#1}{#2}
}
\providecommand{\href}[2]{#2}

%\bibliographystyle{amsplain}
%\bibliography{koszulRef}
\end{document}